\documentclass{amsart}

\usepackage[matrix,arrow]{xy}
\usepackage{amssymb,amsfonts,amsmath,footnote}
\newtheorem{theo}{Th\'eor\`eme}[subsection]
\newtheorem{coro}{Corollaire}
\newtheorem{prop}{Proposition}
\usepackage[francais]{babel}
\begin{document}
\date{Novembre 2015}
\title{Groupes de Galois motiviques et p\'eriodes}
\author{Yves ANDR\'E}
\address{Institut de Math\'ematiques de Jussieu (IMJ-PRG)\\
UPMC\\
UMR 7586 du CNRS\\
\'Equipe de Th\'eorie des Nombres\\
4, place Jussieu\\
Case 247\\
F--75252 Paris Cedex 5}
\email{yves.andre@imj-prg.fr}

 \begin{sloppypar}
\maketitle

\noindent{\bf INTRODUCTION}

\subsection{P\'eriodes.}
  Les int\'egrales ab\'eliennes sont des int\'egrales de 1-formes diff\'erentielles rationnelles sur une courbe alg\'ebrique. Elles d\'ependent du chemin choisi entre les extr\'emit\'es; pour les formes sans p\^ole, les ambigu\"{\i}t\'es sont appel\'ees {\it p\'eriodes ab\'eliennes} car ce sont les (composantes des) {p\'eriodes} des fonctions ab\'eliennes attach\'ees \`a la courbe, par inversion d'Abel-Jacobi. 

Bien que l'interpr\'etation de telles int\'egrales comme p\'eriodes de fonctions 
 fasse d\'efaut 
en dimension sup\'erieure, le terme de {\it p\'eriode} a fini, sur le mode synecdotique, par d\'esigner en g\'eom\'etrie alg\'ebrique toute int\'egrale $\,\int_\Delta \omega\,$ d'une $n$-forme alg\'ebrique $\omega$ prise sur un domaine $\Delta$ limit\'e par des \'equations alg\'ebriques (pour $n=1$, on retrouve les int\'egrales ab\'eliennes).  
 Dans une tradition qui remonte \`a Euler, Legendre et Gauss, deux cas particuliers ont pris une importance consid\'erable: 

- le {\og cas arithm\'etique\fg}, formalis\'e dans \cite{K-Z}: il s'agit de nombres complexes dont les parties r\'eelle et imaginaire sont de la forme $\int_\Delta \omega$, o\`u $\omega$ est une forme diff\'erentielle rationnelle sur une vari\'et\'e alg\'ebrique $X$ d\'efinie sur le corps $\mathbb Q$ des nombres rationnels, et o\`u $\Delta \subset X(\mathbb R)$ est d\'efini par des in\'egalit\'es polynomiales \`a coefficients dans $\mathbb Q$. 

- le {\og cas fonctionnel\fg}: il s'agit de p\'eriodes de formes diff\'erentielles d\'ependant alg\'ebriquement d'un param\`etre $t$ ou de plusieurs, le corps des constantes \'etant $\mathbb C$. Ces fonctions {\og multiformes\fg} sont holonomes \`a croissance mod\'er\'ee (solutions d'\'equations diff\'erentielles lin\'eaires \`a coefficients dans $\mathbb C(t)$ \`a singularit\'es r\'eguli\`eres), fait qui g\'en\'eralise le lien d\'ecouvert par Gauss entre  moyenne arithm\'etico-g\'eom\'etrique et \'equation diff\'erentielle hyperg\'eom\'etrique. 
 
 \subsection{Relations de p\'eriodes.} Dans l'expression $\int_\Delta \omega$ d'une p\'eriode, seul le signe $\int$ n'est pas de nature alg\'ebrique, et l'on s'attend de fait \`a ce que les p\'eriodes soient en g\'en\'eral, dans l'un ou l'autre cas, des nombres ou des fonctions transcendant(e)s. 

On peut s'interroger sur les exceptions\footnote{\`a l'instar de Leibniz, correspondant avec Huygens
 \`a propos du lemme XXVIII des Principia de Newton sur les aires de secteurs d'ovales (cf. e.g. \cite{Va}\cite{Wu}). En langage moderne, ce lemme affirme qu'une telle aire (qui est une p\'eriode dans le cas d'un ovale alg\'ebrique) n'est pas alg\'ebrique en les param\`etres des droites d\'ecoupant le secteur. Newton l'applique au cas de la trajectoire elliptique d'une plan\`ete et en d\'eduit, via la seconde loi de Kepler, que sa position ne d\'epend pas de mani\`ere alg\'ebrique du temps.},
  et plus g\'en\'eralement sur la nature des relations polynomiales entre p\'eriodes (\`a coefficients dans $\mathbb Q$ dans le cas arithm\'etique, dans $\mathbb C(t)$ dans le cas fonctionnel). Les transformations d'une expression $\int_\Delta \omega$ qu'on obtient en jouant avec les propri\'et\'es formelles de $\int\,$,  - \`a savoir la bilin\'earit\'e en $(\Delta, \omega)$, la multiplicativit\'e (Fubini), le changement de variable alg\'ebrique $  \;\int_\Delta f^\ast \omega = \int_{f_\ast \Delta} \omega ,\;$ et la formule de Stokes $ \; \int_\Delta d\omega = \int_{\partial\Delta} \omega \;$ - ,  donnent tautologiquement lieu \`a des relations polynomiales entre p\'eriodes. 
 
 Une r\'eponse tr\`es leibnizienne \`a la question serait alors le philosoph\`eme: {\it\og les propri\'et\'es formelles de $\int$ constituent la raison suffisante des relations polynomiales entre p\'eriodes\fg}.  
  Cette r\'eponse a d'abord \'et\'e propos\'ee par M. Kontsevich, dans le cas arithm\'etique, sous forme d'une conjecture pr\'ecise  
   \cite{K1}, illustr\'ee de maint exemple en collaboration avec D. Zagier \cite{K-Z}. 
  
  Tout r\'ecemment, J. Ayoub a montr\'e que la m\^eme r\'eponse vaut dans le cas fonctionnel, cette fois sous forme d'un th\'eor\`eme, que voici. Notons
     
 - $\mathcal O(\bar{\mathbb D}^n)$ l'alg\`ebre des fonctions analytiques au voisinage du polydisque unit\'e ferm\'e $\,\vert z_i\vert \leq 1\, (i= 1,\ldots, n)$,
 
  - $ \mathcal O_{alg}^\dagger(\bar{\mathbb D}^n)$ le sous-espace de $\mathcal O(\bar{\mathbb D}^n)((t))$ form\'e des s\'eries 
 $    \sum_{i>>-\infty}  h_i(z_1, \ldots , z_n) t^i$ qui sont alg\'ebriques sur $\mathbb C(z_1, \ldots ,z_n, t)$,  et 
  $\mathcal O_{alg}^\dagger(\bar{\mathbb D}^\infty)= \cup\, \mathcal O_{alg}^\dagger(\bar{\mathbb D}^n)$.

\begin{theo}[\cite{Ay5} th. 1.8]  Le noyau de l'application $\mathbb C$-lin\'eaire $$\displaystyle   \sum   h_i  \,t^i\in  \mathcal O_{alg}^\dagger(\bar{\mathbb D}^\infty) \; \;\mapsto \;\;\sum  (\int_{[0,1]^\infty} h_i ) \,  t^i\in \mathbb C((t))$$ est engendr\'e par les \'el\'ements de la forme 
$$ \frac{\partial g}{\partial z_i} - g_{\mid z_i = 1} + g_{\mid z_i = 0} \;\;\;{\it{et}}\;\;\; (f -\int_{[0,1]^\infty} f )\,h$$ o\`u $i\in \mathbb N\setminus 0,\; f, g,h\in O_{alg}^\dagger(\bar{\mathbb D}^\infty)$, $f$ ne d\'epend pas de $t$, et $f$ et $h$ ne d\'ependent pas simultan\'ement d'une m\^eme variable. 
\end{theo}

   \subsection{Groupes de Galois motiviques.} Ceci repose sur de remarquables progr\`es en {\it th\'eorie de Galois motivique}\footnote{la puissance de la th\'eorie motivique dans les questions de relations de p\'eriodes a d\'ej\`a \'et\'e mise en \'evidence dans les travaux de F. Brown \cite{Bro} sur les nombres polyz\^eta (qui sont des p\'eriodes), r\'ecemment recens\'es dans ce s\'eminaire \cite{D4}.}, d\^us principalement \`a M. Nori et surtout J. Ayoub. 
 
 Soit $k$ un sous-corps de $\mathbb C$.  
   Selon A. Grothendieck, il devrait exister une cat\'egorie ab\'elienne $\mathbb Q$-lin\'eaire de motifs mixtes sur $k$, not\'ee ${\bf MM}(k) $, pr\'esentant les trois caract\'eristiques suivantes:
 
 1) \^etre le r\'eceptacle d'une {\og cohomologie universelle\fg}; d'o\`u un foncteur {\og r\'ealisation de Betti\fg} $H_B: {\bf MM}(k) \to {\rm{Vec}}_{\mathbb Q}$ vers les espaces vectoriels de dimension finie sur $\mathbb Q$, qui factorise la cohomologie de Betti,  
  
 2) \^etre {\og de nature g\'eom\'etrique\fg}: ses morphismes devraient \^etre construits {\og en partant de\fg} correspondances alg\'ebriques,
 
 3) \^etre dot\'ee d'un produit tensoriel compatible au produit des vari\'et\'es, qui en fasse une cat\'egorie tannakienne: $H_B$ induirait une $\otimes$-\'equivalence entre ${\bf MM}(k) $ et la cat\'egorie des repr\'esentations de dimension finie d'un $\mathbb Q$-sch\'ema en groupe affine, le groupe de Galois motivique  ${\bf G}_{mot}(k) := {\bf Aut}^\otimes H_B$.
  
 Apr\`es un demi-si\`ecle, on peut consid\'erer ce programme comme d\'esormais accompli - m\^eme si la th\'eorie des motifs est loin d'\^etre achev\'ee\footnote{notamment en ce qui concerne l'interpr\'etation de {\og en partant de\fg} dans le point 2); nous pr\'eciserons de quelles interpr\'etations il s'agit dans la suite, sans sp\'eculer sur des interpr\'etations plus \'etroites.}. Sur la base de nouvelles th\'eories tannakiennes, Nori \cite{N} a construit une cat\'egorie ${\bf MM}_N(k)$ v\'erifiant 1) et 3), et Ayoub \cite{Ay3} une cat\'egorie ${\bf MM}_{Ay}(k)$ v\'erifiant 2)\footnote{avec les pr\'ecautions de la note pr\'ec\'edente.} et 3); en outre, le point 1) permet de construire un $\otimes$-foncteur de la premi\`ere vers la seconde, qui s'av\`ere \^etre une \'equivalence \cite{C-GAS}.  Enfin, d'apr\`es D. Arapura \cite{Ar}, la sous-cat\'egorie form\'ee des objets semi-simples s'identifie \`a la cat\'egorie tannakienne ${\bf M}(k)$ des motifs purs construite ant\'erieurement par le rapporteur en termes de {\og correspondances motiv\'ees\fg}   
 sur les vari\'et\'es projectives lisses \cite{An2}  (elles-m\^emes d\'efinies \`a partir de correspondances alg\'ebriques).

  \subsection{Torseurs des p\'eriodes.}  
  
   Dans \cite{G}, Grothendieck construit la {\it cohomologie de De Rham alg\'ebrique} $H_{dR}^\ast(X)$ d'une $k$-vari\'et\'e alg\'ebrique lisse $X$ comme hypercohomologie du complexe de De Rham alg\'ebrique $\Omega^\ast_X$, et explique que les p\'eriodes proviennent d'un accouplement parfait (apr\`es extension des scalaires \`a $\mathbb C$): 
   $$\; H_{dR}^\ast(X)\otimes H^B_\ast(X)\to \mathbb C,\; (\omega, \Delta)\mapsto  \int_\Delta \omega.\;$$
   Ceci se g\'en\'eralise (non trivialement) au cas non lisse, et m\^eme au cas d'une paire $(X,Y)$, $Y$ \'etant une sous-vari\'et\'e ferm\'ee de $X$, cf. \cite[II 3]{H-MS2}. Par le point 1) ci-dessus, on obtient donc un foncteur {\og r\'ealisation de De Rham\fg} $H_{dR}: {\bf MM}(k) \to {\rm{Vec}}_{k}$, et par le point 3), un ${\bf G}_{mot}(k)_k$-torseur  ${\bf P}_{mot}(k) := {\bf Iso}^\otimes (H_{dR}, H_B\otimes k)$ dot\'e d'un point complexe canonique $\varpi :\; {\rm{Spec}}\, \mathbb C \to {\bf P}_{mot}(k)$ correspondant \`a l'accouplement de p\'eriodes. 
   
   Cette construction montre que le {\it degr\'e de transcendance} sur $k$ des p\'eriodes d'une $k$-vari\'et\'e $X$ est toujours major\'e par la dimension de son groupe de Galois motivique (c'est-\`a-dire l'image de ${\bf G}_{mot}(k)$ dans $GL(H_B(X))$).
 
Lorsque $k = \mathbb Q$ (le cas arithm\'etique), Grothendieck a conjectur\'e que l'image de $\varpi$ est un point g\'en\'erique\footnote{bien qu'il ne l'ait pas publi\'ee, tout comme ses autres id\'ees sur les motifs, une trace de l'existence de cette conjecture (dans le cas des int\'egrales ab\'eliennes) figure en note de bas de page de sa lettre \cite{G}. La premi\`ere version publi\'ee de la conjecture (dans le cas pur) se trouve dans un appendice du livre de S. Lang sur les nombres transcendants \cite{La}. Le rapporteur a pu se rendre compte, lors d'une visite \`a Montpellier durant laquelle J. Malgoire a eu l'obligeance de lui montrer quelques manuscrits de Grothendieck sur les motifs, que cette conjecture y figure bien sous la forme \'enonc\'ee ci-dessus.}, de sorte que la majoration du degr\'e de transcendance des p\'eriodes serait optimale. Compte tenu du point $2)$, c'est une mise en forme du philosoph\`eme: {\it\og toute relation polynomiale entre p\'eriodes est d'origine g\'eom\'etrique\fg}. La compatibilit\'e de ce philosoph\`eme avec le pr\'ec\'edent ne va pas de soi, mais une description pr\'ecise du torseur des p\'eriodes montre que les conjectures de Grothendieck et de Kontsevich sont en fait \'equivalentes (cf. 3.4). Elles forment le socle d'une th\'eorie de Galois des p\'eriodes qui \'etendrait partiellement aux nombres transcendants la th\'eorie de Galois usuelle des extensions alg\'ebriques de $\mathbb Q$ \cite{An4}.
    
\smallskip  Ayoub a propos\'e et {\it d\'emontr\'e} un analogue fonctionnel de la conjecture des p\'eriodes de Grothendieck. En gros, il s'agit de remplacer le groupe de Galois motivique absolu par un groupe relatif qui d\'ecrit les motifs mixtes sur $\mathbb C(t)$ {\og modulo\fg} les motifs d\'efinis sur le sous-corps $\mathbb C$ des constantes. Plus pr\'ecis\'ement, apr\`es r\'eduction des constantes \`a un corps alg\'ebriquement clos $k$ convenable dont le plongement complexe s'\'etend \`a $k(t)$, le groupe qui contr\^ole la $k(t)$-alg\`ebre des p\'eriodes est le noyau ${\bf G}_{mot}(k(t)\mid k)$ de l'homomorphisme canonique ${\bf G}_{mot}(k(t))\to {\bf G}_{mot}(k)$. Les p\'eriodes \'etant holonomes \`a croissance mod\'er\'ee, les relations polynomiales \`a coefficients dans $k(t)$ qui les lient sont contr\^ol\'ees par leur monodromie; le point cl\'e est donc l'existence d'un homomorphisme d'image dense de la monodromie vers ${\bf G}_{mot}(k(t)\mid k)$  \cite{Ay3}, ce qui traduit une version motivique du classique th\'eor\`eme de la partie fixe en th\'eorie de Hodge, cf 5.2.  
  
  Si ce {\og th\'eor\`eme-ma\^{\i}tre\fg} donne en principe toute l'information sur les relations entre p\'eriodes dans le cas fonctionnel en termes de la g\'eom\'etrie sous-jacente, sa traduction concr\`ete dans les 
  situations tr\`es diverses o\`u l'on rencontre
    ces probl\`emes requiert des techniques adapt\'ees, dont nous pr\'esentons un \'echantillon en 1.4 et 5.2.

\section{Relations entre int\'egrales ab\'eliennes (cas fonctionnel).}
Pour entrer en mati\`ere, nous examinons deux probl\`emes d'ind\'ependance d'int\'egrales ab\'eliennes d\'ependant de param\`etres. Ce sera l'occasion de discuter le r\^ole qu'a jou\'e la th\'eorie de Hodge dans ces questions, avant d'indiquer comment celle-ci peut d\'esormais \^etre remplac\'ee, de mani\`ere plus naturelle, par la th\'eorie motivique.  
  
 \subsection{Int\'egrales ab\'eliennes et $1$-motifs.} Consid\'erons une vari\'et\'e ab\'elienne $A$ de dimension $g$, d\'efinie sur un sous-corps $k$ de $\mathbb C$. Son groupe de cohomologie de De Rham $H^1_{dR}(A)$ est de dimension $2g$, engendr\'e par $g$ formes {\og de premi\`ere esp\`ece \fg} $\omega_1, \ldots, \omega_g\in \Omega^1(A)$ et par $g$ formes {\og de seconde esp\`ece \fg} $\eta_1, \ldots, \eta_g$. Si $\gamma_1,\ldots, \gamma_{2g}$ d\'esigne une base du groupe d'homologie $H_1(A(\mathbb C), \mathbb Q)$, la matrice des p\'eriodes de $A$ dans ces bases est la matrice carr\'ee inversible d'ordre $2g$ 
   $\Omega  =  \begin{pmatrix}   
    \int_{\gamma_i} \omega_j & \int_{\gamma_i} \eta_j  
     \end{pmatrix}_{i= 1,\ldots 2g,\, j=1,\ldots g}.$
 
  Tout endomorphisme $f$ de $A$ fournit des relations lin\'eaires entre les coefficients de $\Omega$, du type $\int_{\gamma } f^\ast \omega = \int_{f_\ast\gamma }  \omega$.   Toute polarisation de $A$ donne lieu \`a des relations quadratiques, du type de celles d\'efinissant le groupe des similitudes symplectiques. Il peut exister d'autres relations d'origine g\'eom\'etrique (cf.  2.4). 
  
 Si $x$ est un $k$-point de $A$, on peut aussi consid\'erer les int\'egrales ab\'eliennes\footnote{dans le langage classique, le nom de p\'eriodes est r\'eserv\'e aux $\int_{\gamma_i} \omega_j$, les $\int_{\gamma_i} \eta_j $ sont appel\'es quasi-p\'eriodes, les $\int_0^x  \omega_j$  logarithmes ab\'eliens, tandis que les $\int_0^x  \eta_j$ n'ont pas de nom.} $\int_0^x  \omega_j, \int_0^x \eta_j$, qui d\'ependent du choix d'un chemin entre $0$ et $x$.  Dans le langage des $1$-motifs de P. Deligne et de leur r\'ealisation de Betti et de De Rham \cite{D2}, la matrice augment\'ee (d'une ligne et d'une colonne)
   ${\bf \Omega} = \begin{pmatrix}  1 &  \begin{matrix}  \int_0^x  \omega_j &  \int_0^x \eta_j   \end{matrix} \\
    0 &  \begin{matrix}   
    \int_{\gamma_i} \omega_j & \int_{\gamma_i} \eta_j  
     \end{matrix} \end{pmatrix}_{i= 1,\ldots 2g,\, j=1,\ldots g} $
 est une matrice de p\'eriodes du $1$-motif ${\bf A} = [\mathbb Z \stackrel{1\mapsto x}{\to} A]$. 
  
 Il est loisible de travailler ici avec la cat\'egorie {\it ab\'elienne} des vari\'et\'es ab\'eliennes (resp. $1$-motifs) {\og \`a isog\'enie pr\`es\fg}, obtenue en tensorisant les morphismes avec $\mathbb Q$.

  \subsection{Le th\'eor\`eme du noyau de Manin.} 
 
   Dans la situation relative (le {\og cas fonctionnel\fg}),  consid\'erons un sch\'ema ab\'elien $A\to S$ sur une vari\'et\'e alg\'ebrique complexe lisse, et, en pr\'esence d'une section $x: S\to A$, le $1$-motif ${\bf A} = [\mathbb Z \stackrel{1\mapsto x}{\to} A]$ sur $S$. La cohomologie de De Rham $H^1_{dR}({\bf A}/S)$ (resp. $H^1_{dR}({\bf A}/S)$) est un $\mathcal O_S$-module localement libre muni d'une connexion int\'egrable $\nabla$, dont une premi\`ere construction alg\'ebrique remonte \`a Yu. Manin \cite{M}\footnote{que les int\'egrales ab\'eliennes d\'ependant alg\'ebriquement d'un param\`etre soient solutions d'\'equations diff\'erentielles \`a coefficients polynomiaux en ce param\`etre (comme dans l'exemple hyperg\'eom\'etrique de Gauss) \'etait un fait connu depuis longtemps: on parlait d'\'equations diff\'erentielles de Picard-Fuchs, avant que Grothendieck n'introduise dans \cite{G} la terminologie {\og connexion de Gauss-Manin\fg} - construite peu apr\`es en toute g\'en\'eralit\'e par N. Katz et T. Oda.}; en fait, $(H^1_{dR}({\bf A}/S), \nabla)$ est extension de $(H^1_{dR}({ A}/S),\nabla)$ par $(\mathcal O_S, d)$. On obtient un homomorphisme:
    $$\mu: A(S) \to {\rm{Ext}}^1((H^1_{dR}({ A}/S),\nabla), (\mathcal O_S, d)), \; x\mapsto  (H^1_{dR}({\bf A}/S), \nabla) .$$
Quitte \`a remplacer $S$ par un ouvert dense, on peut supposer $H^1_{dR}({\bf A}/S)$ libre, de base $\omega_1, \ldots, \omega_g\in \Omega^1(A/S),\; \eta_1, \ldots, \eta_g,$ et former comme ci-dessus les matrices de p\'eriodes $\Omega$ (resp. $\bf \Omega$), qui sont des {\og solutions\fg} de  $(H^1_{dR}({ A}/S),\nabla)$  (resp. $(H^1_{dR}({\bf A}/S), \nabla)$) \`a valeurs analytiques multiformes sur $S(\mathbb C)$.  Il est donc \'equivalent de dire que $x$ est dans le noyau de $\mu$, ou bien que le vecteur-ligne de composantes $(\int_0^x \omega_1, \ldots \int_0^x \omega_g,  \int_0^x \eta_1, \ldots \int_0^x \eta_g)$ est combinaison lin\'eaire \`a coefficients constants des vecteurs  $(\int_{\gamma_i} \omega_1, \ldots \int_{\gamma_i} \omega_g,  \int_{\gamma_i} \eta_1, \ldots \int_{\gamma_i} \eta_g), i= 1,\ldots 2g$. 
 C'est le cas, banalement, si $x$ est la section nulle, ou plus g\'en\'eralement (puisque $\mu$ est additive) si $x$ est de torsion. 
 Le th\'eor\`eme du noyau de Manin concerne la r\'eciproque:

\begin{theo} Supposons que $A/S$ n'ait pas de partie fixe (i.e. de sous-sch\'ema ab\'elien non trivial devenant constant sur un rev\^etement \'etale fini de $S$). Alors les conditions suivantes sont \'equivalentes:
\label{cascade} \item $i)$ $x$ est une section de torsion,
\item $ii)$ $(\int_0^x \omega_1, \ldots \int_0^x \omega_g,  \int_0^x \eta_1, \ldots \int_0^x \eta_g)$ est combinaison lin\'eaire \`a coefficients constants des vecteurs  $(\int_{\gamma_i} \omega_1, \ldots \int_{\gamma_i} \omega_g,  \int_{\gamma_i} \eta_1, \ldots \int_{\gamma_i} \eta_g), i= 1,\ldots 2g$,
\item $iii)$ $(\int_0^x \omega_1, \ldots \int_0^x \omega_g )$ est combinaison lin\'eaire \`a coefficients constants des vecteurs  $(\int_{\gamma_i} \omega_1, \ldots \int_{\gamma_i} \omega_g), i= 1,\ldots 2g$.
\end{theo} 
 C'est la forme $iii)$ qui est utilis\'ee par Manin dans sa preuve de la conjecture de Mordell pour les corps de fonctions. Elle se ram\`ene imm\'ediatement \`a $ii)$ lorsque $H^1_{dR}({ A}/S)$ est engendr\'e, en tant que module diff\'erentiel, par les $\omega_i$. Cette condition n'est pas toujours v\'erifi\'ee; or la preuve de \cite{M} semble la supposer implicitement, comme l'a remarqu\'e R. Coleman,  et certaines tentatives ult\'erieures ont achopp\'e sur ce point - voir le compte-rendu de D. Bertrand \cite{Be} (c'est une variante de sa solution que nous pr\'esenterons).

\subsection{Le th\'eor\`eme d'ind\'ependance d'int\'egrales ab\'eliennes modulo les p\'eriodes.} D'apr\`es ce qui pr\'ec\`ede, on a une tour d'extensions diff\'erentielles $\mathbb C(S) \subset \mathbb C(S)(\int_{\gamma_i} \omega_j,\int_{\gamma_i} \eta_j) \subset \mathbb C(S)(\int_{\gamma_i} \omega_j,\int_{\gamma_i} \eta_j)(\int_0^x \omega_j, \int_0^x \eta_j)$.
 Tout comme le th\'eor\`eme du noyau de Manin, le th\'eor\`eme suivant \cite[th.2]{An1} concerne la seconde inclusion.

\begin{theo} Supposons que $A/S$ n'ait pas de partie fixe, et que la section $x$ soit g\'en\'erique (i.e. d'image non contenue dans un sous-sch\'ema en groupe distinct de $A$). Alors les $\int_0^x \omega_j$ et $ \int_0^x \eta_j$ sont alg\'ebriquement ind\'ependants sur $\mathbb C(S)(\int_{\gamma_i} \omega_j,\int_{\gamma_i} \eta_j)$.
\end{theo} 

Par la th\'eorie de Galois diff\'erentielle, le degr\'e de transcendance de l'extension $\mathbb C(S)(\int_{\gamma_i} \omega_j,\int_{\gamma_i} \eta_j)(\int_0^x \omega_j, \int_0^x \eta_j)$ sur $\mathbb C(S)(\int_{\gamma_i} \omega_j,\int_{\gamma_i} \eta_j)$  est la dimension du noyau de l'homomorphisme du groupe de Galois diff\'erentiel de l'extension $ \mathbb C(S)(\int_{\gamma_i} \omega_j,\int_{\gamma_i} \eta_j)(\int_0^x \omega_j, \int_0^x \eta_j)/\mathbb C(S)$ vers celui de $\mathbb C(S)(\int_{\gamma_i} \omega_j,\int_{\gamma_i} \eta_j) /\mathbb C(S)$. Comme les connexions en jeu ont des singularit\'es r\'eguli\`eres \`a l'infini, c'est aussi, par la correspondance de Riemann-Hilbert et le choix d'un point base $s\in S(\mathbb C)$, la dimension du noyau $U_{\bf A} = {\rm{Ker}}\,({\bf H}\to H)$ de l'homomorphisme (surjectif) entre groupes de monodromie point\'es en $s$.  Par d\'efinition, le groupe de monodromie $H$ (resp. $\bf H$) est l'adh\'erence de Zariski de l'image de $\pi_1(S(\mathbb C), s)$ dans $GL(H_1(A_s))$ (resp. $GL(H_1({\bf A}_s))$).

\subsection{R\^ole des variations de structures de Hodge.} Pour contr\^oler ces groupes, il faut prendre en compte l'origine g\'eom\'etrique de ces repr\'esentations de monodromie. On retiendra que ces syst\`emes locaux sont sous-jacents \`a des {\it variations de structures de Hodge rationnelles polarisables} (pure pour $A/S$, mixte pour ${\bf A}/S$). Cette propri\'et\'e suffira car la cat\'egorie des sch\'emas ab\'eliens (resp. des $1$-motifs) sur $S$ \`a isog\'enie pr\`es est une sous-cat\'egorie ab\'elienne {\it pleine} de celle ${\bf VSHP}_S$ (resp. ${\bf VSHMP}_S$) des variations de structures de Hodge rationnelles (resp. mixtes) polarisables sur $S$  \cite[4.4.3]{D1}\cite[10.1.3]{D2}.
 La cat\'egorie ${\bf VSHMP}_S$ est tannakienne, neutralis\'ee par le foncteur fibre en $s$. On peut donc associer aux  VSHMP  attach\'ees \`a $A$ et $\bf A$ des groupes tannakiens $G$ et $\bf G$, sous-groupes alg\'ebriques de $GL(H_1(A_s))$ et de  
 $GL(H_1({\bf A}_s))$ qui contiennent $H$ et $\bf H$ respectivement. On a homomorphisme surjectif ${\bf G} \to G$ qui prolonge ${\bf H}\to H$.

Un r\'esultat central de la th\'eorie des VSHP est le th\'eor\`eme de la partie fixe (d\^u \`a Deligne dans le {\og cas g\'eom\'etrique\fg} et \`a W. Schmid en g\'en\'eral, cf. \cite[th. 11]{P-S}): 
 {\it le plus grand sous-syst\`eme local constant du syst\`eme local sous-jacent \`a une} VSHP {\it est sous-jacent \`a une sous}-VSHP. Autrement dit, pour tout ${ \underline V}\in {\bf VSHP}_S$, on a une structure de Hodge induite sur ${ \underline V}_s^{\pi_1(S(\mathbb C), s)}$, ind\'ependante de $s$.  Ce r\'esultat s'\'etend aux VSHMP  {\og admissibles\fg}, celles v\'erifiant une certaine condition \`a l'infini \cite{S-Z}, qui est satisfaite dans le cas des $1$-motifs \cite[lem. 5]{An1}.  
   On en tire les cons\'equences suivantes: 
 
 $a)$ toute VSHP (ou VSHMP admissible) dont le syst\`eme local sous-jacent est constant est elle-m\^eme constante. En particulier, {\it tout sch\'ema ab\'elien $A\to S$ dont la monodromie est triviale est constant}; 
  
 $b)$ le groupe de monodromie de toute VSHMP admissible est un sous-groupe normal de son groupe tannakien (cf. \cite[lem. 1]{An1}); les repr\'esentations de ce groupe tannakien qui se factorisent par la monodromie correspondent \`a des VSHMP constantes\footnote{de sorte que le foncteur fibre en un point g\'en\'eral $s$ identifie en fait le groupe quotient au groupe de Mumford-Tate de cette fibre. Le th. 5.1 ci-dessous en est l'analogue motivique (plus profond).}. En particulier $H$ (resp. $\bf H$) {\it est un sous-groupe normal de} $G$ (resp. $\bf G$).
 
 $c)$ $H_1(A_s)^{H}$ est stable sous $G$, donc {\it correspond \`a un sous-sch\'ema ab\'elien constant de} $A/S$ (\`a isog\'enie pr\`es), nul si $A/S$ est sans partie fixe;
  
 $d)$   Soit $\bf V$  une $\bf G$-repr\'esentation extension de la repr\'esentation triviale $\mathbb Q$ par une repr\'esentation $V$ du quotient $G$ telle que $V^H=0$. Alors $\bf V$ est scind\'ee si et seulement si elle l'est en tant que $\bf H$-repr\'esentation: en effet, $\mathbb Q$ se rel\`eve en la $\bf G$-repr\'esentation ${\bf V}^{\bf H}$, n\'ecessairement triviale. Ceci s'applique \`a $(V = H_1(A_s), {\bf V} = H_1({\bf A}_s))$ si $A/S$ est sans partie fixe.

\subsection{Preuve des th\'eor\`emes 1.1 et 1.2 (esquisse).} 

L'implication $i)\Rightarrow iii)$ du th\'eor\`eme 1.1 est banale. Prouvons $ii)\Rightarrow i)$: via le plongement des 1-motifs sur $S$ \`a isog\'enie pr\`es dans ${\bf VSHMP}_S$, il s'agit de voir que l'extension 
de VSHMP $0\to H_1(A/S) \to H_1({\bf A}/S)\to \mathbb Q\to 0$ est scind\'ee si et seulement si celle des syst\`emes locaux sous-jacents l'est (du moins si $A/S$ est sans partie fixe), ce qui, en traduction tannakienne, est le point $d)$ ci-dessus.  

\smallskip Passons au th\'eor\`eme 1.2. Quitte \`a remplacer $S$ par un rev\^etement \'etale fini, on peut supposer que l'anneau End$_S A$ est l'anneau des endomorphismes $R $ d'une fibre g\'en\'erale $A_s$; et que le point $x_s $ n'est contenu dans aucun sous-sch\'ema en groupe propre de $A_s$.  Ceci implique que pour tout $f\in R\setminus 0$, $f(x_s)$ n'est pas de torsion. 

 Reprenons les notations du point $d)$ ci-dessus. Soit $m\in {\bf V}$ un point se projetant sur l'\'el\'ement $1$ du quotient $\mathbb Q$ de ${\bf V}$. La r\`egle $u\mapsto  u(m) - m $ d\'efinit un plongement $\bf G$-\'equivariant de $U_{\bf A}(\mathbb Q)$ dans $V$, ind\'ependant du choix de $m$ (construction de Bashmakov-Ribet, cf. \cite[prop. 1]{An1}). Son image est un sous-$\bf G$-module $W\subset V$; il correspond \`a un sous-sch\'ema ab\'elien de $A/S$ \`a isog\'enie pr\`es. Si $V\neq W$, sa fibre en $s$ est contenue dans le noyau d'un \'el\'ement $f\in R\setminus 0$.  Le $1$-motif \`a isog\'enie pr\`es ${{\bf A}^f} = [\mathbb Z \stackrel{1\mapsto f(x )}{\to} A ]$ est dans la cat\'egorie ab\'elienne engendr\'ee par $\bf A $, 
 donc $U_{{\bf A}^f}$ est quotient de $U_{\bf A}= {\rm{Ker}}({\bf H}\to H)$, et en fait nul par construction de $f$.  Ceci contredit l'implication $ii)\Rightarrow i)$ du th\'eor\`eme 1.1  appliqu\'e \`a la section $f(x)$, donc $W=V$ et par suite ${\rm{dim}}\,U_{\bf A} = {\rm{dim}}_{\mathbb Q}\,V = 2g$. 

\smallskip  Pour terminer, prouvons l'implication $iii)\Rightarrow ii)$ du th\'eor\`eme 1.1. Soit $\mathcal M$ le sous-module diff\'erentiel de $H^1_{dR}(A/S)$ engendr\'e par les $\omega_j$. Il s'agit de montrer que si  $\mu(x) $ est non nul, sa restriction $\mu(x)_{\mathcal M}$ \`a ${\rm{Ext}}^1(\mathcal M, (\mathcal O_S, d))$ l'est aussi. 
 Quitte \`a remplacer $x$ par un multiple, on peut supposer que le plus petit sous-sch\'ema en groupe de $A$ contenant l'image de $x$ est un sous-sch\'ema ab\'elien, et quitte \`a remplacer $A$ par ce dernier, on peut supposer $x$ g\'en\'erique. Il suffit alors de voir que pour tout sous-module diff\'erentiel non nul $\mathcal M$ de $H^1_{dR}(A/S)$, $\mu(x)_{\mathcal M}\neq 0$.  Or $\mathcal M$ et l'extension par $(\mathcal O_S, d))$ restriction de $H^1_{dR}({\bf A}/S)$ correspondent alors \`a des quotients $V_{\mathcal M}$ et ${\bf V}_{\mathcal M}$ de $V_{\mathbb C}$ et ${\bf V}_{\mathbb C}$ respectivement. Notons $U_{\mathcal M}$ l'image de $U_{{\bf A}, \mathbb C}$ dans $GL({\bf V}_{\mathcal M})$. L'isomorphisme $U_{\bf A}(\mathbb C) \stackrel{\sim}{\to} V_{\mathbb C}$ obtenu dans la d\'emonstration du th\'eor\`eme 1.2  compos\'e avec la projection $V_{\mathbb C}\to V_{\mathcal M}$ se factorise \`a travers $U_{\mathcal M}(\mathbb C)$. Si $\mathcal M$ est non nul, il en est de m\^eme de $V_{\mathcal M}$ donc aussi de $U_{\mathcal M}$, de sorte que $\mu(x)_{\mathcal M}\neq 0$. \qed

\smallskip {\it Remarque}.  L'intervention un peu artificielle des structures de Hodge, efficace parce que vari\'et\'es ab\'eliennes et $1$-motifs (\`a isog\'enie pr\`es) se laissent d\'ecrire en ces termes, peut \^etre \'evit\'ee en les rempla\c cant par des motifs. On sait en effet plonger les $1$-motifs dans une cat\'egorie tannakienne de motifs (\S 3.3), et on dispose d'analogues motiviques du th\'eor\`eme de la partie fixe (\S 5.1) permettant de remplacer dans les arguments ci-dessus les groupes $G, \bf G$ par des groupes de Galois motiviques. On peut alors r\'einterpr\'eter ces arguments comme l'\'elaboration concr\`ete, dans certaines situations particuli\`eres de relations entre int\'egrales ab\'eliennes, du th\'eor\`eme d'Ayoub mentionn\'e en 0.4, voire envisager de les \'etendre hors du cadre ab\'elien. Plus qu'un point technique, il s'agit d'un changement de {\og bon cadre\fg} pour \'etudier les relations entre p\'eriodes.

\section{Groupes de Galois motiviques purs.}
\subsection{La construction de Grothendieck.}

  Dans toute la suite de cet expos\'e, $k$ {\it d\'esigne un sous-corps de} $\mathbb C$ et $\bar k$ sa fermeture alg\'ebrique. Les motifs purs sur $k$ sont ceux attach\'es aux $k$-vari\'et\'es projectives lisses; on s'attend \`a ce qu'ils forment une cat\'egorie ab\'elienne {\it semi-simple}, refl\'etant les propri\'et\'es cohomologiques de puret\'e/semi-simplicit\'e, connues ou attendues, de ces vari\'et\'es.  
 
 La construction conjecturale bien connue de Grothendieck repose sur l'id\'ee que le monde des motifs purs est celui de la g\'eom\'etrie \'enum\'erative (combinatoire de configurations en g\'eom\'etrie projective). Il consid\`ere la cat\'egorie mono\"{\i}dale ayant pour objets les $k$-vari\'et\'es projectives lisses, et pour morphismes les correspondances alg\'ebriques (\`a coefficients rationnels et de dimension \'egale \`a celle du but),  modulo l'\'equivalence num\'erique. Passant \`a l'enveloppe karoubienne, la classe de ${\bf P}^1$ se d\'ecompose en ${\mathbb Q} \oplus {\mathbb Q}(-1)$, et en $\otimes$-inversant ${\mathbb Q}(-1)$, on obtient  une cat\'egorie mono\"{\i}dale ${\bf M}_{G}(k)$ dont tout objet est dualisable. 
Grothendieck montre alors, sous ses {\og conjectures standard\fg}, que ${\bf M}_{G}(k)$ est ab\'elienne semi-simple, que les cohomologies se factorisent \`a travers des $\otimes$-foncteurs de source ${\bf M}_{G}(k)$, et - en introduisant \`a cette occasion les {\it cat\'egories tannakiennes} - que ${\bf M}_{G}(k)$ est \'equivalente \`a la cat\'egorie des repr\'esentations de dimension finie d'un $\mathbb Q$-sch\'ema en groupe affine ${\bf G}_{mot}^{pur}(k)$, le groupe de Galois motivique pur.

L'hypoth\`eque originelle des conjectures standard\footnote{les conjectures standard
 apparaissent pour la premi\`ere fois dans une lettre de Grothendieck \`a J.-P. Serre de 1965 \cite{G-S} qui se termine par ces mots: {\og Ce qu'il faut pour le moment, c'est inventer un proc\'ed\'e pour d\'eformer un cycle de dimension pas trop grande, pour le pousser \`a l'infini. Peut-\^etre aurais-tu envie d'y r\'efl\'echir de ton c\^ot\'e? Je viens seulement de m'y mettre aujourd'hui-m\^eme, et t'\'ecris faute de trouver une id\'ee.\fg}. Cinquante ans plus tard, on en est au m\^eme point. Heureusement, ces questions fondamentales {\it ne sont pas au fondement} de la th\'eorie des motifs, comme le montre son \'evolution ult\'erieure (si les conjectures standard \'etaient fausses, aucun \'enonc\'e de cet expos\'e n'en serait affect\'e!).}a notoirement nui \`a l'essor et \`a la r\'eputation de la th\'eorie, qui n'est sortie des limbes qu'au d\'ebut des ann\'ees 80 sous l'impulsion de Deligne, en sacrifiant provisoirement son essence g\'eom\'etrique pour se  concentrer sur les {\og syst\`emes de r\'ealisations\fg}.
 Ainsi,  l'\'etude des p\'eriodes incite \`a consid\'erer la cat\'egorie tannakienne ${\rm{Vec}}_{k,\mathbb Q}$ dont les objets sont des triplets $(W\in {\rm{Vec}}_k, \,V\in {\rm{Vec}}_{\mathbb Q}, \,\iota: W\otimes_k \mathbb C \stackrel{\sim}{\to} V\otimes_{\mathbb Q} \mathbb C)$ et la sous-cat\'egorie ab\'elienne engendr\'ee par ceux de la forme $(H_{dR}(X), H_B(X), \varpi)$ o\`u $X$ est une $k$-vari\'et\'e projective lisse. 
  La th\'eorie de Hodge absolue met en jeu, quant \`a elle, la cat\'egorie dont les objets sont des familles $(W\in {\rm{Vec}}_k, \,V_\sigma \in {\bf{SH}}, \,\iota_\sigma: W\otimes_k \mathbb C \stackrel{\sim}{\to} V_\sigma \otimes_{\mathbb Q} \mathbb C)_{\sigma: k \to \mathbb C}$, o\`u ${\bf{SH}}$ d\'esigne la cat\'egorie des structures de Hodge rationnelles, etc...  Il ne s'agit pas \`a proprement parler de cat\'egories de motifs: les morphismes, d\'efinis purement en termes d'alg\`ebre lin\'eaire, ne sont pas de nature g\'eom\'etrique; pour autant, leur \'etude a contribu\'e \`a fa\c conner le formalisme de Galois motivique et men\'e \`a de brillantes applications \cite{D3}\cite{D-M-O-S}. 
  
    Le retour aux motifs proprement dits a eu lieu avec le th\'eor\`eme de U. Jannsen \cite{J}
     affirmant, ind\'ependamment des conjectures standard, que ${\bf M}_G(k)$ est ab\'elienne semi-simple... ce qui ne suffit toutefois pas \`a construire les r\'ealisations ni  ${\bf G}_{mot}^{pur}(k)$.

\subsection{Correspondances motiv\'ees et motifs purs.} La modification minimale de la construction de Grothendieck conduisant \`a une th\'eorie de Galois motivique pure inconditionnelle (et conservant la nature g\'eom\'etrique des morphismes) est expos\'ee dans \cite{An2}. Sur $k$, les conjectures standard se r\'eduisent \`a celle qui pr\'edit que, pour une vari\'et\'e polaris\'ee $X$ de dimension $d$, l'inverse de l'isomorphisme de Lefschetz $ H^i(X) \stackrel{l_X}{\to} H^{2d-i}(X)$ est donn\'e, comme $l_X$ lui-m\^eme (cup-produit it\'er\'e avec la polarisation), par une correspondance alg\'ebrique. En introduisant les {\it correspondances motiv\'ees}, obtenues \`a partir des correspondances alg\'ebriques modulo \'equivalence homologique en {\og adjoignant \fg} formellement les $l_X^{-1}$ (la notion ne d\'epend pas des polarisations choisies), et en suivant la construction de Grothendieck, on obtient une {\it cat\'egorie $\mathbb Q$-lin\'eaire tannakienne semi-simple de motifs purs} ${\bf M}(k)$, \`a travers laquelle se factorisent les cohomologies de Weil. D'o\`u un foncteur fibre {\og r\'ealisation de Betti pure\fg} $H_B:  {\bf M}(k) \to {\rm{Vec}}_{\mathbb Q}$ et un $\mathbb Q $-sch\'ema en groupe affine ${\bf G}_{mot}^{pur}(k):= {\bf Aut}^\otimes H_B$, de sorte que $H_B$ s'enrichit en une  $\otimes$-\'equivalence ${\bf M}(k)\stackrel{\sim}{\to} {\rm{Rep}}\,{\bf G}_{mot}^{pur}(k)$. Tout ceci repose sur le th\'eor\`eme de l'indice de Hodge (et se transpose au cas relatif \cite{Ar-D}).

 \`A toute $k$-vari\'et\'e projective lisse $X$ (et plus g\'en\'eralement tout motif pur), on associe ainsi un ${\mathbb Q}$-groupe alg\'ebrique r\'eductif, son {\it groupe de Galois motivique} ${\bf G}(X)\subset GL(H_B(X))$, qui n'est autre que l'image de ${\bf G}_{mot}^{pur}(k)$. Cette construction permet en principe de ramener les probl\`emes de nature {\og motivique\fg} \`a des questions de th\'eorie des repr\'esentations des groupes r\'eductifs\footnote{de m\^eme que la th\'eorie de Galois usuelle ram\`ene les probl\`emes d'extensions finies de $k$ \`a des questions d'actions de groupes finis; au reste, Gal$(\bar k/k)$ est un quotient de ${\bf G}_{mot}^{pur}(k)$, correspondant aux vari\'et\'es de dimension $0$.}.

 Le th\'eor\`eme de la partie fixe de Deligne a un avatar motivique (qui s'en d\'eduit): {\it si $X\to S$ est un morphisme projectif lisse, et $s\in S(\mathbb C),\;H_B(X_s)^{\pi_1(S(\mathbb C), s)}$ est la r\'ealisation de Betti d'un sous-motif du motif de $X$}; de mani\`ere \'equivalente,  $H_B(X_s)^{\pi_1(S(\mathbb C), s)}$ est stable sous ${\bf G}(X_s)$ \cite{An2}.

 \subsection{La conjecture des p\'eriodes pures.} La cohomologie de De Rham alg\'ebrique $H_{dR}(X) := {\bf H}_{Zar}(\Omega^\ast_X)$ fournit une r\'ealisation $H_{dR}:  {\bf M}(k) \to {\rm{Vec}}_k$. \`A toute $k$-vari\'et\'e projective lisse $X$ (et plus g\'en\'eralement \`a tout motif pur), on associe un torseur ${\bf P}(X):= {\bf Iso}^\otimes(H_{dR}, H_{B}\otimes k)_{<X>_\otimes}$ sous  ${\bf G}(X)_k$, muni d'un point complexe canonique $\varpi_X: {\rm{Spec}}\,\mathbb C \to {\bf P}(X)$ donn\'e par l'isomorphisme de p\'eriodes. 
 
 La conjecture des p\'eriodes de Grothendieck dans le cas pur\footnote{dans une formulation ind\'ependante des conjectures standard. On prendra garde par ailleurs \`a bien distinguer la conjecture des p\'eriodes, qui a trait \`a la r\'ealisation des p\'eriodes ${\bf M}(k) \to {\rm{Vec}}_{k, \mathbb Q}$ (et implique sa pleine fid\'elit\'e si $k\subset \bar{\mathbb Q}$), de la conjecture de Hodge, qui a trait \`a la r\'ealisation de Hodge ${\bf M}(k) \to {\bf SH}$ (et implique sa pleine fid\'elit\'e si $k= \bar k$). Elles sont de nature enti\`erement diff\'erente, cf. \cite[ch. 7]{An3}.}  pr\'edit que  {\it si $k\subset \bar{\mathbb Q}$},
  {\it l'image de $\varpi_X$ est un point Zariski-dense de $ {\bf P}(X)$},
 
 ou, de mani\`ere \'equivalente, 
 
 {\it $ {\bf P}(X)$ est connexe et le degr\'e de transcendance des p\'eriodes de $X$ est}  dim ${\bf G}(X)$.
 
 Par exemple, ${\bf G}({\bf P}^1)=\mathbb G_m$ et les p\'eriodes sont $1$ et $2\pi i$, de sorte que la conjecture \'equivaut dans ce cas \`a la transcendance de $\pi$. Pour un panorama des r\'esultats de la th\'eorie des nombres transcendants en faveur de cette conjecture, voir \cite{Wal}\footnote{pour une g\'en\'eralisation de la conjecture au cas d'un corps de base $k$ non n\'ecessairement alg\'ebrique, qui implique entre autres la conjecture de Schanuel sur les exponentielles, voir \cite[23.4.1]{An3}: le degr\'e de transcendance sur $\mathbb Q$ de la $k$-alg\`ebre engendr\'ee par les p\'eriodes de $X$ serait toujours minor\'ee par la dimension de ${\bf G}(X)$.}.  
  
\subsection{Cas des p\'eriodes ab\'eliennes.} Le seul r\'esultat g\'en\'eral dans cette direction est le th\'eor\`eme de G. W\"ustholz (version ab\'elienne du th\'eor\`eme de Baker sur les logarithmes):

\begin{theo}[cf. \cite{Wu}] Si $k = \bar{\mathbb Q}$, toute relation $k$-{\emph{lin\'eaire}} entre p\'eriodes de $1$-formes sur une $k$-vari\'et\'e ab\'elienne provient de ses endomorphismes.    
\end{theo}  

Pour d'autres cas de la conjecture qui s'en d\'eduisent (indirectement), voir \cite{B-C}. 

\smallskip Pour aller plus loin dans l'analyse de la nature motivique des relations polynomiales entre p\'eriodes ab\'eliennes, commen\c cons par pr\'eciser la nature des motifs ab\'eliens. En partant du th\'eor\`eme de la partie fixe motivique et en adaptant des arguments de Deligne \cite{D-M-O-S}, on peut d\'emontrer une variante affaiblie de la conjecture de Hodge pour les vari\'et\'es ab\'eliennes:
 
 \begin{theo}[\cite{An2}]  Si $k=\bar k$, la r\'ealisation de Hodge fait de la cat\'egorie tannakienne des motifs purs engendr\'ee par les vari\'et\'es ab\'eliennes une sous-cat\'egorie pleine de ${\bf SH}$. En particulier, tout cycle de Hodge sur une $k$-vari\'et\'e ab\'elienne est motiv\'e\footnote{il y a ici une subtilit\'e: la conjecture standard est connue pour les $k$-vari\'et\'es ab\'eliennes (Grothendieck-Lieberman), mais il ne s'ensuit pas que tout cycle motiv\'e soit alg\'ebrique sur ces vari\'et\'es; le point est que les cycles motiv\'es font intervenir $l_Y^{-1}$ pour des vari\'et\'es auxiliaires $Y$ non n\'ecessairement ab\'eliennes - dans la situation du th\'eor\`eme, il s'agit de pinceaux ab\'eliens compacts.}. 
 \end{theo}  

Pour toute $k$-vari\'et\'e ab\'elienne $X$, ${\bf G}(X)$ co\"{\i}ncide donc avec le groupe de Mumford-Tate de $X$ (groupe tannakien attach\'e \`a la structure de Hodge $H_1(X)$). Si $X$ est \`a multiplication complexe (ce qui permet de supposer $k = \bar{\mathbb Q}$), il s'agit d'un tore dont le groupe de caract\`eres se calcule par une recette explicite en termes du type CM. Par ailleurs, dans le cas de multiplication complexe par un corps cyclotomique (ou plus g\'en\'eralement une extension ab\'elienne de $\mathbb Q$),  les p\'eriodes s'expriment comme produits de valeurs de la fonction $\Gamma$ d'Euler en des nombres rationnels, ce qui g\'en\'eralise la formule classique de Lerch-Chowla-Selberg.

\begin{prop}[\cite{An3}, 24.6] Restreinte aux vari\'et\'es ab\'eliennes \`a multiplication complexe cyclotomique, la conjecture des p\'eriodes de Grothendieck \'equivaut \`a la {\emph{conjecture de Lang-Rohrlich}}: toute relation polynomiale \`a coefficients dans $\mathbb Q$ entre valeurs de $\Gamma$ en des rationnels provient des \'equations fonctionnelles de $\Gamma$. 
\end{prop} 
 
 Du fait que la fonction $\Gamma$ n'est pas elle-m\^eme de nature motivique (elle n'est li\'ee aux motifs ab\'eliens qu'\`a travers ses valeurs aux points rationnels), cette {\og traduction\fg} de la conjecture des p\'eriodes est tr\`es indirecte: reposant sur le th. 2.2, elle s'adosse aussi  \`a une longue s\'erie de travaux sur les p\'eriodes ab\'eliennes de type CM (G. Anderson, Deligne \cite{D3}, B. Gross, G. Shimura, A. Weil \cite{We}...) ainsi que sur l'analyse des relations de distributions, telles celles v\'erifi\'ees par $\Gamma$ (D. Kubert).

Cette {\og traduction\fg} peut d'ailleurs \^etre lue dans l'autre sens: c'est au fond la nature motivique m\^eme des relations entre valeurs de $\Gamma$ que l'on exploite pour calculer rapidement sur machine des relations telles que 
$\, \frac{\Gamma(\frac{1}{24})   \Gamma(\frac{11}{24}) }{\Gamma(\frac{5}{24})   \Gamma(\frac{7}{24}) } =   \sqrt{6+3\sqrt{3}}$, 
 cf. \cite{B-Z}.

\section{Groupes de Galois motiviques mixtes selon Nori.} 

\subsection{Th\'eorie tannakienne de Nori.}  Soit $F$ un corps de caract\'eristique nulle (corps de coefficients, qu'on peut supposer \'egal \`a $\mathbb Q$ dans toute la suite).  La th\'eorie tannakienne invent\'ee par Grothendieck en vue de la construction de groupe de Galois motivique, et d\'evelopp\'ee par N. Saavedra, part d'une cat\'egorie $F$-lin\'eaire mono\"{\i}dale sym\'etrique $\mathcal T$ (avec End ${\bf 1}= F$) dont tout objet est dualisable (i.e. la tensorisation avec un objet quelconque a un adjoint), et munie d'un foncteur $F$-lin\'eaire mono\"{\i}dal $f: \mathcal T \to {\rm{Vec}}_F$ vers les $F$-espaces vectoriels de dimension finie. Elle donne un crit\`ere pour que $f$ s'enrichisse en une \'equivalence mono\"{\i}dale $\mathcal T \stackrel{\sim}{\to} {\rm{Rep}}_F G$ pour un $F$-sch\'ema en groupe affine $G$:  $\mathcal T$ doit \^etre ab\'elienne et $f$ fid\`ele exact.  

Comme l'ont montr\'e Nori puis Ayoub, on peut affaiblir consid\'erablement ces conditions si l'on se contente d'enrichir $f$  en un foncteur universel $\mathcal T \to {\rm{Rep}}_F G$, mais pas n\'ecessairement une \'equivalence. 

Dans la construction de Nori, on peut m\^eme partir d'un carquois (i.e. graphe orient\'e) quelconque $\mathcal Q$, plut\^ot que d'une cat\'egorie. Une  repr\'esentation de $\mathcal Q$ \`a valeurs dans une cat\'egorie ab\'elienne $\mathcal A$ associe \`a tout sommet de $\mathcal Q$ un objet de $\mathcal A$  et \`a toute fl\`eche un morphisme entre les objets associ\'es \`a sa source et \`a son but respectivement. Les repr\'esentations de  $ \mathcal Q$ dans $\mathcal A$ forment de mani\`ere \'evidente une cat\'egorie.

L'\'enonc\'e de base est le suivant \cite{N} (d\'evelopp\'e dans  \cite{Ar}\cite{Bru}\cite{H-MS2}\cite{BV-C-L}):

\begin{prop} A toute $F$-repr\'esentation $f: \mathcal Q \to {\rm{Vec}}_F$, associons la $F$-cog\`ebre $\mathcal C = \varinjlim_{\mathcal Q' \,{\rm{fini}}\, \subset \mathcal Q}\, ({\rm{End}} f_{\mid \mathcal Q'})^\vee$. Alors $f$ s'enrichit en une repr\'esentation $  \mathcal Q \to {\rm{Comod}}_F \mathcal C$ \`a valeurs dans la cat\'egorie ab\'elienne des $\mathcal C$-comodules de dimension finie. Elle est universelle parmi les enrichissements de $f$ en une repr\'esentation $  \mathcal Q \to \mathcal A$  vers une cat\'egorie ab\'elienne $F$-lin\'eaire, tels que le foncteur d'oubli $\mathcal A \to {\rm{Vec}}_F$ soit exact et fid\`ele. 
\end{prop} 

(La colimite filtrante sur les sous-carquois finis pallie le d\'efaut de dualit\'e entre alg\`ebres et cog\`ebres en dimension infinie). 
  La th\'eorie explicite ensuite les conditions mono\"{\i}dales sur $f$ qui font de $\mathcal C$ une big\`ebre, voire une alg\`ebre de Hopf (variation subtile sur le th\`eme: tout sous-mono\"{\i}de Zariski-ferm\'e de $GL_n$ est un groupe, cf. \cite[7.3.6]{H-MS2}).

\subsection{Motifs mixtes de Nori.} Gr\^ace \`a la proposition 3.1, Nori construit sa cat\'egorie ab\'elienne des motifs mixtes comme le r\'eceptacle d'une cohomologie universelle, par le biais d'un carquois qui code les propri\'et\'es standard de toute cohomologie relative. Ses sommets sont des triplets $(X, Y, i)$ form\'es d'une $k$-vari\'et\'e $X$, d'une sous-vari\'et\'e ferm\'ee $Y$, et d'un entier $i\in \mathbb N$. Ses fl\`eches sont de deux types: 
 celles $(X, Y , i) \to (X',Y',i)$ provenant de morphismes de paires, et celles $(X, Y, i) \to (Y, Z, i-1)$ qu'on associe aux triplets de sous-vari\'et\'es embo\^{\i}t\'ees $Z\subset Y\subset X$. La repr\'esentation $f$ est celle induite par la cohomologie relative $H^i_B(X,Y, \mathbb Q) = H^i(X(\mathbb C), j_! \mathbb Q)$ (o\`u $j$ d\'esigne l'inclusion de $X(\mathbb C)\setminus Y(\mathbb C)$ dans $X(\mathbb C)$). Par 3.1. on obtient ainsi une cog\`ebre $\mathcal H_N^{eff}$ sur $F= \mathbb Q$. 
 
Pour en faire une big\`ebre, la formule de K\"unneth 
  $$H^j_B(X \times X', Y\times X' \cup Y'\times X, \mathbb Q) = \bigoplus_{i+i'= j} H^i_B(X,Y, \mathbb Q) \otimes H^{i'}_B(X',Y', \mathbb Q) $$ incite \`a poser $(X, Y, i)\otimes  (X', Y', i') := (X\times X', Y\times X' \cup Y'\times X, i+i')$. Toutefois, du fait de la sommation, on n'a pas compatibilit\'e de $f$ \`a $\otimes$ ... sauf dans le {\it cas cellulaire}, c'est-\`a-dire quand les paires $(X, Y)$ et $(X', Y')$ n'ont qu'un seul groupe de cohomologie non nul. La cl\'e est l'existence d'un analogue alg\'ebrique de la filtration par le squelette d'un complexe simplicial, construite en it\'erant la proposition suivante\footnote{{\og Nori's basic lemma\fg}, d\'emontr\'e ind\'ependamment, dans un cadre plus g\'en\'eral, par A. Beilinson et par K. Vilonen.}:
  
  \begin{prop}[\cite{N}, \cite{H-MS2} 2.5] Soient $X$ une $k$-vari\'et\'e affine de dimension $n$, et $Z$ un ferm\'e de dimension $<n$. Il existe un ferm\'e $Y$ de dimension $<n$ contenant $Z$ tel que $H^i_B(X,Y, \mathbb Q) = 0$ si $i\neq n$.
  \end{prop}

\noindent{\sc Preuve} --- On peut supposer $X\setminus Z$ lisse. Par r\'esolution des singularit\'es, on peut trouver $\tilde X$ projective lisse, $D\subset X$ diviseur \`a croisements normaux et $\pi: \tilde X \setminus D\to X$ propre surjectif, et un isomorphisme au-dessus de $X\setminus Z$;  on peut supposer en outre que $\tilde Z$ de $Z$ soit un diviseur \`a croisement normaux coupant $D$ transversalement.  Pour une section hyperplane g\'en\'erale $H$ de $\tilde X$, posons $D' = H \cup \tilde Z$. Alors $D\cup D' $ est encore un diviseur \`a croisements normaux et $\tilde X\setminus D'$ est affine. 

 Posons $Y= \pi(D' \setminus D\cap D')$. C'est un ferm\'e de la vari\'et\'e affine $X$, de dimension $<n$ et contenant $Z$.  D'apr\`es le lemme d'annulation d'Artin (cf. \cite[2.3.8]{H-MS2}), $H^i_B(X,Y, \mathbb Q)=0$ si $i>n$. Supposons $i<n$. Par excision (cf. \cite[2.1.7]{H-MS2}), on a $H^i_B(X,Y, \mathbb Q)= H^i_B(\tilde X \setminus D, D'\setminus (D\cap D'), \mathbb Q)$, que la dualit\'e de Poincar\'e pour les paires  met en dualit\'e avec $H^{2n-i}_B(\tilde X \setminus D', D\setminus (D\cap D'), \mathbb Q)$ (cf. \cite[2.4.5]{H-MS2}), qui est nul d'apr\`es Artin. \qed

\smallskip
Par un argument de suite spectrale, on peut en d\'eduire que ${\rm{Comod}}\,\mathcal H_N^{eff}$ est engendr\'ee par les images par $f$ de paires cellulaires, puis (en faisant attention \`a changer le signe dans la sym\'etrie selon la r\`egle de Koszul), que ${\mathcal H}_N^{eff}$ est une big\`ebre \cite[ch. 3]{H-MS2}. 
 Pour obtenir une alg\`ebre de Hopf $\mathcal H_N$, il faut en outre inverser $\mathbb Q(-1)$, l'image de $(\mathbb P^1, \emptyset, 2)$ (ce qui se fait plus commod\'ement en modifiant le carquois pour int\'egrer les torsions de Tate). Le $\mathbb Q$-sch\'ema en groupe d\'efini par son spectre est le {\it groupe de Galois motivique de Nori} ${\bf G}_{mot}^N(k)$, et $${\bf MM}_{N}(k) := {\rm{Rep}}_{\mathbb Q}\, {\bf G}_{mot}^{N}(k) = {\rm{Comod}}\,\mathcal H_{N}$$ la cat\'egorie tannakienne des {\it motifs mixtes de Nori} sur $k$ \`a coefficients dans $\mathbb Q$\footnote{tout ceci s'\'etend aux coefficients entiers \cite{H-MS2}, et au cas relatif \cite{Ar}.}.  

  \smallskip  On renvoie \`a  \cite[th. 8.1.9]{H-MS2} pour la formulation pr\'ecise de l'universalit\'e de ${\bf MM}_N(k) $.     

\subsection{Motifs purs et $1$-motifs revisit\'es.} Les morphismes de ${\bf MM}_N(k)$ \'etant d\'efinis formellement \`a partir de propri\'et\'es abstraites des cohomologies relatives, il n'est pas du tout clair qu'ils soient {\og de nature g\'eom\'etrique\fg}. Voici deux r\'eponses positives partielles (une autre, plus compl\`ete, sera donn\'ee au \S 4.5):

\begin{theo}[\cite{Ar} 6.4] La cat\'egorie ${\bf M}(k)$ des motifs purs (\S 2.2) s'identifie canoniquement \`a la sous-cat\'egorie tannakienne de ${\bf MM}_N(k) $ form\'ee des objets semi-simples, et la dualit\'e tannakienne identifie ${\bf G}_{mot}^{pur}(k)$ au plus grand quotient pro-r\'eductif de ${\bf G}_{mot}^{N}(k)$.
\end{theo}

Les morphismes entre motifs de Nori attach\'es \`a des vari\'et\'es projectives lisses sont donc les correspondances motiv\'ees (la preuve commence par la construction d'une filtration par le poids sur ${\bf MM}_N(k) $). 

\begin{theo}[\cite{Ay-BV}] La sous-cat\'egorie ab\'elienne de ${\bf{MM}}_N(k)$ engendr\'ee par les images des $(X,Y,i)$ avec $i\leq 1$ s'identifie \`a la cat\'egorie des $1$-motifs de Deligne.
\end{theo}

 \subsection{Torseurs des p\'eriodes selon Nori-Kontsevich.} La cohomologie de De Rham alg\'ebrique s'\'etend en une cohomologie relative, pour des paires $(X, Y)$ non n\'ecessairement lisses, cf. \cite[II 3]{H-MS2}; apr\`es tensorisation par $\mathbb C$, elle devient canoniquement isomorphe \`a la cohomologie de Betti. Elle donne donc lieu \`a un $\otimes$-foncteur  $H_{dR}:\,{\bf{MM}}_N(k)\to Vec_k$ (r\'ealisation de De Rham), et l'on peut derechef construire le torseur ${\bf P}_{mot}^N(k) := {\bf Iso}^\otimes(H_{dR}, H_{B}\otimes k)$ sous ${\bf G}_{mot}^{N}(k)_k$, dot\'e de son point canonique $\varpi :\; {\rm{Spec}}\, \mathbb C \to {\bf P}_{mot}^N(k)$.

  Les fonctions sur ce torseur s'interpr\`etent comme {\og p\'eriodes abstraites\fg}\footnote{r\'esultat annonc\'e par Kontsevich \cite{K1}\cite{K2} - qui l'attribue \`a Nori - , et d\'emontr\'e en d\'etail dans \cite[III]{H-MS2}; voir aussi \cite{K-Z}, \cite{H-MS1}. On trouvera par ailleurs dans \cite[III]{H-MS2} une comparaison d\'etaill\'ee des variantes de notions de p\'eriodes qu'on rencontre dans la litt\'erature; voir aussi \cite{BB}.}. 
  Suivant Kontsevich-Zagier,  on consid\`ere le $k$-espace vectoriel $\mathcal P^{eff}_{KZ}$ engendr\'e par les symboles $[X, Y, i,  \gamma, \omega]$, o\`u  $( X, Y, i)$ est un sommet du carquois de Nori, $\gamma \in H_i(X(\mathbb C),Y(\mathbb C),\mathbb Q)$ et $ \omega\in H^i_{dR}(X,Y )$, modulo les relations suivantes:
  
  - $(\mathbb Q, k)$-lin\'earit\'e en $(\gamma, \omega)$, 
  
  - (changement de base) pour $f: X'\to X$ tel que $f(Y')\subset Y , \; [X, Y, i, f_\ast \gamma', \omega] = [X', Y', i,  \gamma', f^\ast\omega] $,

  - (Stokes) pour $Z\subset Y \subset X,$ $ \,\gamma \in H_i(X(\mathbb C),Y(\mathbb C),\mathbb Q)$ et $\omega\in H^{i-1}_{dR}(Y,Z )$, $\; [X, Y, i,    \gamma, d\omega] = [Y, Z, i-1,  \partial\gamma,  \omega] $. 
  
 C'est en fait une $k$-alg\`ebre, et posant $\underline{2\pi i}:= [\mathbb G_m, \emptyset, 1, S^1, \frac{dt}{t}]$, on d\'efinit la $k$-alg\`ebre des {\it p\'eriodes abstraites}  $\mathcal P_{KZ}$ comme $\mathcal P^{eff}_{KZ}[\underline{2\pi i}^{-1}]$.  
 
  L'int\'egrale  $ \int_\gamma \omega$ ne d\'epend que de la classe de $[X, Y, i, f_\ast \gamma, \omega] $ (ce qui ne fait que traduire les r\`egles du calcul int\'egral) et envoie $\underline{2\pi i}$ sur $2\pi i$, d'o\`u un homomorphisme $k$-lin\'eaire   $\;\int:\; \mathcal P_{KZ} \to \mathbb C.$  
   
 \begin{theo}[\cite{K2}\cite{H-MS2} ch. 12] ${\bf P}_{mot}^{N}(k)$ est canoniquement isomorphe \`a ${\rm{Spec}}\,\mathcal P_{KZ}$, en sorte que $\varpi$ correspond \`a $\int$.
 \end{theo} 
 
  \begin{coro} Supposons $k\subset \bar{\mathbb Q}$. Alors  $\varpi$ est un point g\'en\'erique {\emph{(conjecture de Grothendieck)}} si et seulement si $\int$ est injectif {\emph{(conjecture de Kontsevich-Zagier)}}.
 \end{coro} 
  
  \noindent {\it Remarque.} Les nombres polyz\^eta $ \zeta({\underline s})= \sum_{n_1>\ldots > n_k}\,  n_1^{-s_1}\ldots n_k^{-s_k}$ sont des p\'eriodes (pour $k=\mathbb Q)$. D'apr\`es Brown \cite{Bro}, ce sont en fait les p\'eriodes des motifs de Tate mixtes sur $\mathbb Z$, et ils s'expriment comme combinaison lin\'eaire rationnelle de polyz\^etas de Hoffman (ceux o\`u $s_i = 2$ ou $3$). Si comme on s'y attend, les motifs de Tate mixtes sur $\mathbb Z$ forment une sous-cat\'egorie tannakienne de celle de Nori, la restriction de la conjecture des p\'eriodes de Grothendieck-Kontsevich-Zagier \`a ces motifs se traduirait par l'ind\'ependance lin\'eaire sur $\mathbb Q$ des polyz\^etas de Hoffman.

 \section{Groupes de Galois motiviques mixtes selon Ayoub.}
 
\subsection{Cat\'egories triangul\'ees de motifs. } Le programme motivique de Grothendieck envisageait des motifs sur une base plus g\'en\'erale que Spec $k$ ainsi que l'existence d'un formalisme des 4 op\'erations de changement de base, dans un cadre d\'eriv\'e. 

 Cette partie ambitieuse du programme est elle aussi d\'esormais accomplie. Dans un premier temps (apr\`es des \'etudes pr\'eliminaires sur les syst\`emes d\'eriv\'es de r\'ealisations, puis les travaux pionniers de A. Beilinson, S. Bloch, A. Suslin...), trois versions de la {\it cat\'egorie triangul\'ee des motifs} sur $k$ ont \'et\'e construites, qui se sont av\'er\'ees \'equivalentes (M. Hanamura, M. Levine, V. Voevodsky) \cite{V-S-F}\cite{Le}. Ensuite, sur les traces de Voevodsky et F. Morel, Ayoub a mis en place le formalisme des op\'erations de Grothendieck \cite{Ay1} (d\'evelopp\'e ult\'erieurement aussi dans les travaux de D.-C. Cisinski et F. D\'eglise).

Pour b\^atir sa version de la th\'eorie de Galois motivique, Ayoub se place dans la cat\'egorie triangul\'ee non born\'ee ${\bf{DM}}(k)$ des motifs de Voevodsky sur $k$, \`a coefficients rationnels (la version ${\bf{DM}}(S)$ sur une base lisse $S$ intervient aussi de mani\`ere transitoire). Grosso modo, on part de la cat\'egorie des $k$-vari\'et\'es lisses avec pour morphismes les correspondances finies \`a coefficients dans $\mathbb Q$, et on consid\`ere les foncteurs $\mathbb Q$-lin\'eaires contravariants de cette cat\'egorie vers ${\rm{Vec}}_{\mathbb Q}$ ({\it pr\'efaisceaux avec transferts}).  On obtient la {\it cat\'egorie triangul\'ee des motifs mixtes effectifs} 
comme sous-cat\'egorie de la cat\'egorie d\'eriv\'ee des pr\'efaisceaux avec transferts form\'ee des objets v\'erifiant la descente \'etale\footnote{ou Nisnevich, cela revient au m\^eme \`a coefficients rationnels \`a cause des transferts. Prendre garde qu'il s'agit d'un formalisme covariant: indexation {\og homologique\fg} des complexes; derri\`ere cette convention, il y a le choix de ce qu'on consid\`ere comme motif effectif: $H^i$ ou $H_i$, $\mathbb Q(-1)$ ou $\mathbb Q(1)$?} et l'invariance par $\mathbb A^1$-homotopie (elle se r\'ealise alternativement comme localisation de la cat\'egorie d\'eriv\'ee des faisceaux \'etales avec transferts). Pour obtenir  ${\bf{DM}}(k)$, il faut ensuite $\otimes$-inverser\footnote{il s'agit l\`a d'une inversion non na\"{\i}ve, selon la m\'ethode des spectres comme en th\'eorie de l'homotopie stable.}  le motif de Tate $\mathbb Q(1) := {\mathbb G}_m[-1]$. C'est une cat\'egorie mono\"{\i}dale sym\'etrique, et dans la plus petite sous-cat\'egorie triangul\'ee ${\bf{DM}}_{gm}(k)$ stable par facteurs directs et torsion de Tate $\otimes {\mathbb Q}(n)$ et contenant les $k$-vari\'et\'es lisses, les objets sont dualisables: en fait ${\bf{DM}}_{gm}(k)$ est engendr\'ee par les (facteurs directs des) motifs des $k$-vari\'et\'es projectives lisses.

\subsection{Th\'eorie tannakienne d'Ayoub.} Cette th\'eorie vise, comme celle de Nori, \`a enrichir un $\otimes$-foncteur $f: \mathcal T \to {\rm{Vec}}_F$ en un foncteur $\mathcal T \to {\rm{Rep}}_F G$. Toutefois, comme la condition principale est l'existence d'adjoints, il y a lieu de s'affranchir de la finitude et travailler avec les Ind-cat\'egories ${\rm{VEC}}_F$ (espaces vectoriels de dimension quelconque sur $F$) et ${\rm{REP}}_F G$. 
Un \'enonc\'e typique est le suivant: 

\begin{prop}  Soient $\mathcal T$ une cat\'egorie $F$-lin\'eaire mono\"{\i}dale sym\'etrique et $f: \mathcal T\to {\rm{VEC}}_F$ un $\otimes$-foncteur $F$-lin\'eaire. On suppose que $f$ admet un adjoint \`a droite $g$ commutant aux sommes directes infinies. Alors $fg({\bf 1})$ est canoniquement munie d'une structure de big\`ebre, et $f$ s'enrichit en un $\otimes$-foncteur de $ \mathcal T$ vers les $fg({\bf 1})$-comodules. 
\end{prop} 

 Dans le cas o\`u $\mathcal T= {\rm{REP}}_F G$ et $f$ le foncteur d'oubli, $g$ est la tensorisation avec $\mathcal O(G)$, et on retrouve $G=$ Spec $ fg({\bf 1})$.  
 Dans un cadre \'elargi, l'\'enonc\'e pr\'ecis utilis\'e pour construire les groupes de Galois motiviques est le suivant:
 
 \begin{prop}[\cite{Ay3} 1.5]  Soit $f: \mathcal T\to \mathcal V$ un $\otimes$-foncteur entre cat\'egories mono\"{\i}dales sym\'etriques. Supposons que $f$ admette un adjoint \`a droite $g$, ainsi qu'une section mono\"{\i}dale $e$ qui admet elle-m\^eme un adjoint \`a droite $u$. 
 Supposons encore que 

$(\ast)$  pour tout  $T\in \mathcal T$ et tout $V\in \mathcal V$, le morphisme canonique compos\'e 
 $c_{V,T}: $
 
 \noindent $g(V) \otimes T \to gf(g(V)\otimes T) \cong g(fg(V)\otimes f(T))\to g(V\otimes f(T))$ soit un isomorphisme. 

Alors $fg({\bf 1})$ est canoniquement muni d'une structure d'alg\`ebre de Hopf dans $\mathcal V$, et $f$ s'enrichit en un $\otimes$-foncteur de $ \mathcal T$ vers les $fg({\bf 1})$-comodules.  
 \end{prop} 

Nous nous contenterons d'\'ecrire les formules pour la comultiplication et la structure de $fg({\bf 1})$-comodule sur $f(T)$ ($u$ sert pour l'antipode):

- $ fg({\bf 1}) \to fgfg({\bf 1}) = fg({\bf 1} \otimes fefg {\bf 1}) \stackrel{c^{-1}_{{\bf 1}, efg {\bf 1}}}{\to} f(g{\bf 1} \otimes efg {\bf 1}) =  fg{\bf 1} \otimes fefg {\bf 1}  = fg{\bf 1} \otimes  fg{\bf 1}  $,  

 - $ f(T) \to fgf(T) = fg({\bf 1} \otimes fef(T))\stackrel{c^{-1}_{{\bf 1}, ef T}}{\to}  f(g{\bf 1} \otimes ef(T)) = fg{\bf 1} \otimes f(T) $.

 \subsection{Motifs mixtes d'Ayoub.} Ayoub applique 4.2. en prenant pour $f$ un {\og foncteur de Betti\fg} convenable $B^\ast: {\bf{DM}}(k) \to D({\rm{VEC}}_{\mathbb Q})$ ($k$ \'etant un corps plong\'e dans  $\mathbb C$).
 Pour le construire, il consid\`ere l'analogue analytique complexe ${\bf{DM}}^{an}$ de ${\bf{DM}}(k)$, o\`u les $k$-vari\'et\'es lisses sont remplac\'ees par des vari\'et\'es analytiques complexes, $\mathbb A^1$ par le disque unit\'e ouvert $\mathbb D^1$, et la topologie \'etale par la topologie usuelle.  Tirant parti de la simplicit\'e des hyperrecouvrements d'une vari\'et\'e analytique par des polydisques $\mathbb D^n$, il montre que ${\bf{DM}}^{an}$ est canoniquement \'equivalente \`a  $D({\rm{VEC}}_{\mathbb Q})$. Le foncteur compos\'e $B^{\ast}$ envoie le motif d'une $k$-vari\'et\'e lisse sur le complexe de cha\^{\i}nes singuli\`eres de $X(\mathbb C)$ dans $D({\rm{VEC}}_{\mathbb Q})$; plus g\'en\'eralement, pour un complexe de pr\'efaisceaux avec transferts $\mathcal F^\ast$, $B^\ast(\mathcal F^\ast)$  se calcule comme complexe total associ\'e \`a $\mathcal F^\ast(\bar{\mathbb D}^\ast_{et})$, o\`u $\bar{\mathbb{D}}^*_{et}$ est un certain objet cocubique construit en termes des voisinages \'etales des polydisques ferm\'es $\bar{\mathbb D}^n $ dans $\mathbb A^n_{\mathbb C}\,$ \cite[2.2]{Ay3}. 

Comme le $\otimes$-foncteur d'analytification $ {\bf{DM}}(k)\to {\bf{DM}}^{an}$ provient d'une adjonction de Quillen, $B^\ast$ admet un adjoint \`a droite not\'e $B_\ast$; il admet aussi une section mono\"{\i}dale \'evidente et toutes les conditions de 4.2 s'av\`erent remplies,
  ce qui donne naissance \`a une alg\`ebre de Hopf $\,  B^\ast B_\ast \mathbb Q\,$  dans  $D({\rm{VEC}}_{\mathbb Q})$.   

\begin{prop}[\cite{Ay3} 2.3] $H_i( { B^\ast B_\ast \mathbb Q})= 0 $ pour $i<0$.  
\end{prop} 

\noindent{\sc Indication sur la preuve} --- D'apr\`es \cite{L-W}, le complexe de De Rham donne naissance \`a un objet $\Omega^\ast_{/k}$ de ${\bf{DM}}(k)$\footnote{C'est aussi une cons\'equence d'un r\'esultat profond de Morel-Cisinski-D\'eglise-Ayoub selon lequel ${\bf{DM}}(k)$ est \'equivalente \`a la cat\'egorie analogue {\og sans transferts\fg}, cf. \cite[app. B]{Ay3}; c'est dans ce cadre-ci, plus flexible, que bien des constructions se font, notamment celle des 4 op\'erations de Grothendieck.}. Le th\'eor\`eme de Grothendieck s'interpr\`ete alors comme un isomorphisme entre $\Omega^\ast_{/k}$ et $ B_\ast \mathbb C$ (apr\`es tensorisation par $\mathbb C$). Il s'agit donc de voir la propri\'et\'e d'annulation pour $B^\ast \Omega^\ast_{/k}$, qui se calcule comme ${\rm{Tot}}\,\Omega^\ast(\bar{\mathbb D}^\ast_{et})$. On aboutit  \`a un {\og complexe de De Rham infini\fg}, nul en degr\'es (homologiques) $<0$:

\noindent $  \stackrel{d}{\to} \Omega^{\infty - j}     \stackrel{d}{\to} \Omega^{\infty - {(j-1)}}  \to \cdots \Omega^{\infty - 0} \to 0 $, avec 
 $\displaystyle{\Omega^{\infty - j} := \bigoplus_{I\subset \mathbb N\setminus 0, \vert I\vert = j} \, \mathcal O^{(I)}_{alg}(\bar{\mathbb D}^\infty)\bigwedge_{i\in \mathbb N\setminus (I\cup 0) }  dz_i }$,   o\`u   - 
  $\mathcal O_{alg}(\bar{\mathbb D}^n)$ est l'alg\`ebre des fonctions analytiques sur le polydisque ferm\'e qui sont alg\'ebriques sur $k(z_1,  \ldots, z_n)$,

\noindent  - $\mathcal O_{alg}(\bar{\mathbb D}^\infty) = \bigcup \mathcal O_{alg}(\bar{\mathbb D}^n)$, et 
  
\noindent  - $\mathcal O^{(I)}_{alg}(\bar{\mathbb D}^\infty)$ est le sous-espace des fonctions s'annulant pour $z_i = 0$ ou $1$ si $i\in I$.
 
  \begin{coro} $H_0(B^\ast B_\ast \mathbb Q)  $ h\'erite d'une structure de $\mathbb Q$-alg\`ebre de Hopf $H_{Ay}$. \end{coro} 

Le $\mathbb Q$-sch\'ema en groupe d\'efini par son spectre est le {\it groupe de Galois motivique d'Ayoub} ${\bf G}_{mot}^{Ay}(k)$ , et $${\bf MM}_{Ay}(k) := {\rm{Rep}}_{\mathbb Q}\, {\bf G}_{mot}^{Ay}(k) = {\rm{Comod}}\,\mathcal H_{Ay}$$ est la cat\'egorie tannakienne des {\it motifs mixtes d'Ayoub} sur $k$ \`a coefficients dans $\mathbb Q$.

\subsection{Torseurs des p\'eriodes selon Ayoub.} La cohomologie de De Rham alg\'ebrique donne lieu \`a un $\otimes$-foncteur  ${\bf{MM}}_{Ay}(k)\to Vec_k$ (r\'ealisation de De Rham), et l'on peut de nouveau construire le torseur ${\bf P}_{mot}^{Ay}(k) := {\bf Iso}^\otimes(H_{dR}, H_{B}\otimes k)$ sous ${\bf G}_{mot}^{Ay}(k)_k$, dot\'e de son point canonique $\varpi :\; {\rm{Spec}}\, \mathbb C \to {\bf P}_{mot}^{Ay}(k)$. Son alg\`ebre de fonctions n'est autre que le $0$-i\`eme groupe d'homologie du complexe de De Rham infini ci-dessus. Explicitement, soit $\mathcal P_{Ay}^{eff} $ le quotient de $\mathcal O_{alg}(\bar{\mathbb D}^\infty) = \cup\, \mathcal O_{alg}(\bar{\mathbb D}^n)$ par le sous-$k$-espace engendr\'e par les \'el\'ements de la forme
$$  \frac{ \partial g}{\partial z_i} - g_{\mid z_i =1} + g_{\mid z_i = 0}  \; \; \; (i\in \mathbb N\setminus 0).$$  C'est une $k$-alg\`ebre, et on d\'efinit $\mathcal P_{Ay} $  en inversant la classe d'un \'el\'ement convenable de  $\mathcal O_{alg}(\bar{\mathbb D}^1)$ dont l'int\'egrale sur $[0,1]$ vaut $2\pi i$. 

  Par la r\`egle de Newton-Leibniz, l'int\'egration  $h\mapsto  \int_{{[0,1]}^\infty} \, h  \,dz_1dz_2\cdots $ passe au quotient et d\'efinit un homomorphisme $k$-lin\'eaire   $\,\int_{\square}:\; \mathcal P_{Ay} \to \mathbb C.$   
 
 \begin{theo}[\cite{Ay5}] ${\bf P}_{mot}^{Ay}(k)$ est canoniquement isomorphe \`a ${\rm{Spec}}\,\mathcal P_{Ay}$, en sorte que $\varpi$ correspond \`a $\int_{\square}$.
 \end{theo} 
 
   \begin{coro} Supposons $k\subset \bar{\mathbb Q}$. Alors  $\varpi$ est un point g\'en\'erique {\emph{(conjecture de Grothendieck)}} si et seulement si $\int_\square$ est injectif {\emph{(variante d'Ayoub de la conjecture de Kontsevich-Zagier)}}.  \end{coro} 
  
  \noindent {\it Remarque.} Ce qui surprend quand on compare avec le cadre de Nori-Kontsevich-Zagier, c'est d'une part que les relations de Stokes n'y figurent que sous leur forme la plus \'el\'ementaire (Newton-Leibniz), et surtout que les relations de changement de base ont disparu. En fait, elles se d\'eduisent de Stokes. Ayoub propose de s'en convaincre, dans le cas d'une seule variable $z_1$, par un petit calcul ing\'enieux: partant de $h(z_1)\in \mathcal O_{alg}(\bar{\mathbb D}^1)$ et d'une fonction $f(z_1)$ alg\'ebrique qui envoie le disque unit\'e
   dans lui-m\^eme en fixant $0$ et $1$, la formule de changement de variable montre que  $g(z_1) :=  f'(z_1)h(f(z_1)) - h(z_1)$ est dans le noyau de $\int_{\square}$. On peut l'\'ecrire sous la forme pr\'edite, en passant \`a deux variables: posant $f_1= f(z_1)-z_1,\; f_2= -z_2f'(z_1)+ z_2-1,$ et $g_i = f_i\cdot h(z_2f(z_1)+(1-z_2)z_1)$ pour $i=1,2$, on  obtient $g_{1\mid z_1 = 0 \,{\rm{ou}} \,1} = 0, \, g_{2\mid z_2 = 0} = - h(z_1), \, g_{2\mid z_2 = 1} =  f'(z_1)h(f(z_1))$, et $\frac{ \partial g_1}{\partial z_1} + \frac{ \partial g_2}{\partial z_2}= 0$, d'o\`u  $g(z_1) = \sum_1^2\, (\frac{ \partial g_i}{\partial z_i} - g_{i \mid z_i =1} + g_{i\mid z_i = 0} ) $.

\subsection{Equivalence des cat\'egories de Nori et d'Ayoub.} Annonc\'ee dans \cite{Ay5} et \cite{Ay6}, elle est d\'emontr\'ee en d\'etail dans \cite{C-GAS}.

\begin{theo}[\cite{C-GAS}] On a une $\otimes$-\'equivalence canonique ${\bf MM}_N(k)   \stackrel{\sim}{\to} {\bf MM}_{Ay}(k) $.
\end{theo}

   Le foncteur provient de l'universalit\'e des motifs de Nori. Son inverse est construit via les torseurs de p\'eriodes, en observant qu'on a un homomorphisme canonique $\mathcal P_{Ay}\to \mathcal P_{KZ}$, compatible avec les \'evaluations $\int$ et $\int_\square$: si $A$ est une $k[z_1, \ldots, z_n]$-alg\`ebre \'etale contenant $g\in \mathcal O_{alg}(\bar{\mathbb D}^n)$, et $Y$ est le diviseur de $X := {\rm{Spec}}\, A$ donn\'e par $\prod z_i(z_i-1)= 0$, on envoie la classe de $g$ dans $\mathcal P_{Ay}$ sur celle de $[X, Y, n, \gamma, g dz_1\cdots dz_n]$, o\`u $\gamma$ est donn\'ee par $[0,1]^n \hookrightarrow \bar{\mathbb D}^n \to X(\mathbb C)$. 

\begin{coro} $\mathcal P_{Ay}\cong \mathcal P_{KZ}$.
\end{coro}

\noindent (C'est donc bien de la {\og m\^eme\fg} conjecture de Grothendieck qu'il s'agit dans les cor. 3.6 et 4.6). 
 
  \noindent {\it Remarque.} Le formalisme d'Ayoub permet d'\'ecrire les p\'eriodes $\int_{\Delta} \omega$ comme combinaisons $k$-lin\'eaires de p\'eriodes dont le domaine d'int\'egration est un cube\footnote{c'est bien s\^ur ici que sont cach\'ees les relations de changement de variable alg\'ebrique qui n'apparaissent plus dans le th. 0.1.}. Dans l'autre sens, on peut fixer plut\^ot $\omega$ et exprimer les p\'eriodes comme combinaisons $k$-lin\'eaires de volumes de {\og solides alg\'ebriques\fg} \cite{Y}.

 \smallskip
 Dans la suite, nous identifierons gr\^ace \`a 4.7 les deux cat\'egories de motifs mixtes, et omettrons les indices $N, \,Ay$.
  On peut r\'esumer toute la situation dans le diagramme essentiellement commutatif suivant (qui, insistons, ne d\'epend d'aucune conjecture), sorte d'{\og organigramme\fg} de la th\'eorie actuelle des motifs:
  $$\begin{matrix}   {\bf{CHM}}(k) &  {\hookrightarrow} & {\bf{DM}}_{gm}(k) & \to & D^b({\bf MM}(k) )\\ \downarrow & &&& \downarrow \sum H^i \\ {\bf{M}}(k) & &\hookrightarrow &&  {\bf MM}(k)
       \end{matrix}$$ 
(${\bf{CHM}}(k)$ est la cat\'egorie des {\og motifs de Chow\fg} \`a coefficients rationnels, le foncteur vertical de gauche le passage de l'\'equivalence rationnelle \`a l'\'equivalence homologique, et le second foncteur horizontal du haut est construit par la th\'eorie de Nori).

  \section{Les th\'eor\`emes d'Ayoub sur les p\'eriodes fonctionnelles.}
  
  \subsection{Le {\og th\'eor\`eme de la partie fixe\fg} motivique.} Soit $K$ un corps de fonctions sur $k= \bar k$ et supposons que le plongement complexe de $k$ se prolonge  en un plongement complexe $\sigma: K\hookrightarrow \mathbb C$. Le changement de base $\otimes_k K$ induit un homomorphisme ${\bf G}_{mot}(K)\to {\bf G}_{mot}(k)$ qui est surjectif \cite[2.34]{Ay3}. Notons ${\bf G}_{mot}(K\mid k)$ son noyau. Par ailleurs, soit $ \pi_1((K\mid k)^{an}  )$ le groupe pro-discret limite du pro-syst\`eme des groupes fondamentaux (bas\'es en $\sigma$) des espaces de points complexes des mod\`eles lisses de $K/k$. 
   
  \begin{theo}[\cite{Ay3} 2.57]  On a un homomorphisme canonique  $\pi_1((K\mid k)^{an}  ) \to  {\bf G}_{mot}(K\mid k )$
  dont l'image est Zariski-dense.
  \end{theo} 
  
 \noindent {\it Remarque.} Comme en 1.4 a), on en d\'eduit que si $S$ est une $k$-vari\'et\'e de corps de fonctions $K$, tout syst\`eme local motivique $\underline M\in {\bf{MM}}(S)$, dont le syst\`eme local sous-jacent est constant, est constant (i.e. provient d'un motif sur $k$).   
 
 \noindent{\sc Indication sur la preuve de 5.1} ---  On se ram\`ene au cas o\`u $K/k$ est de degr\'e de transcendance $1$. Le premier pas consiste \`a d\'ecrire $ {\bf G}_{mot}(K\mid k )$ \`a l'aide de la construction tannakienne de la prop. 4.2 . Pour cela, on consid\`ere une $k$-courbe lisse point\'ee $(S, s)$ de corps de fonctions $K$ (et de morphisme structural not\'e $p: S\to {\rm{Spec}}\, k)$). On applique 4.2 au foncteur $s^\ast: DM_{sm}(S)\to DM(k)$ (o\`u $DM_{sm}(S)$ est la plus petite cat\'egorie triangul\'ee de $DM(S)$ contenant les objets dualisables et stable par coproduit), pour obtenir une alg\`ebre de Hopf $\mathcal H(S,s)$. On obtient l'alg\`ebre de Hopf $ {\bf G}_{mot}(K\mid k )$ en appliquant $B^\ast$ et en passant \`a la colimite sur les mod\`eles $S$ (en prenant garde aux changements de point base). Le point cl\'e ici est que $B^\ast \mathcal H(S,s)$ est concentr\'e en degr\'e $0$ (ce qui rend $ {\bf G}_{mot}(K\mid k )$ plus accessible que les groupes absolus $ {\bf G}_{mot}(K)$ ou $ {\bf G}_{mot}( k )$). 
 
 Le deuxi\`eme pas consiste \`a exhiber une application canonique de $H_0(B^\ast \mathcal H(S,s))$ vers l'alg\`ebre des fonctions sur l'enveloppe proalg\'ebrique de $\pi_1(S(\mathbb C), s)$ et \`a montrer qu'elle est injective. L\`a, le point cl\'e est que les sections globales du syst\`eme local $\mathcal L$ sous-jacent \`a un objet $\underline M\in DM_{sm}(S)$ se calculent comme Im$(p^{an\ast}p^{an}_\ast  \mathcal L \to \mathcal L)$, et que l'image de $ p^{\ast}p_\ast  \underline M \to \underline M$ a un sens motivique (des arguments alternatifs, dans le contexte des motifs de Nori-Arapura, ont \'et\'e propos\'es ind\'ependamment par Nori et P. Jossen).

   \subsection{Analogue fonctionnel de la conjecture des p\'eriodes \`a la Grothendieck.}  
   
   Revenons \`a la situation de p\'eriodes d\'ependant alg\'ebriquement d'un param\`etre $t$ comme dans l'introduction. Elles sont associ\'ees \`a une paire $(X, Y)$ au-dessus d'une courbe alg\'ebrique $S$ munie d'une coordonn\'ee $t$, ou plus g\'en\'eralement (quitte \`a restreindre $S$)  \`a un syst\`eme local motivique $\underline M\in {\bf{MM}}(S)$. On peut descendre le corps de base $\mathbb C$ \`a un sous-corps $k= \bar k$ sur lequel $S$ et $\underline M$ sont d\'efinis, et tel que le plongement complexe de $k$ se prolonge en un plongement complexe $\sigma$ du corps de fonctions $K$ de la courbe descendue. Comme on l'a d\'ej\`a vu, les p\'eriodes de $M$ sont des fonctions holonomes multiformes \`a croissance mod\'er\'ee (0.1, 0.4, 1.3), solutions de la connexion de Gauss-Manin associ\'ee \`a $\underline M$. Par la th\'eorie de Galois diff\'erentielle et la correspondance de Riemann-Hilbert, le th. 5.1 implique 
   
   \begin{coro}[\cite{Ay6} th. 44] Le degr\'e de transcendance de l'extension de $\mathbb C(t)$ engendr\'ee par les p\'eriodes de $\underline M$ est \'egal \`a la dimension de l'image de ${\bf G}_{mot}(K\mid k)$ dans $GL(H_B(\underline M_\sigma))$. 
   \end{coro} 
  
{\it Remarque 1.} En pratique, conna\^itre le degr\'e de transcendance ne suffit pas: il faut faire intervenir la g\'eom\'etrie du torseur des p\'eriodes, selon que l'on a affaire \`a des probl\`emes d'ind\'ependance de p\'eriodes fonctionnelles sur $\mathbb C$ plut\^ot que sur $\mathbb C(t)$ (cf 1.2), ou de p\'eriodes relativement \`a d'autres (cf 1.3), ou de certaines composantes seulement de la matrice des p\'eriodes (cf \cite[ch. 23]{An3}). En g\'en\'eral, groupes et torseurs ne suffisent pas pour analyser ces questions, et il faut faire intervenir les {\it vari\'et\'es quasi-homog\`enes} pour contr\^oler les relations polynomiales entre solutions de Gauss-Manin \cite{An5}.

 \subsection{Analogue fonctionnel de la conjecture des p\'eriodes \`a la Kontsevich.}  Dans la m\^eme situation, l'usage d'une r\'ealisation de Betti attach\'ee \`a $\sigma$ comme ci-dessus n'est pas tout \`a fait satisfaisante, du fait de la descente de $\mathbb C$ \`a $k$ et du choix arbitraire de $\sigma$, et surtout de ce que les accouplements $H_{dR}(X_{k(t)})\otimes H_B(X_\sigma)\to \mathbb C$ sont \`a valeurs constantes, donn\'ees par les valeurs des p\'eriodes en $\sigma$. Pour prendre en compte le d\'eveloppement en s\'erie des p\'eriodes, il y a lieu de remplacer $\sigma$ par un {\og point tangentiel\fg} et de construire la {\it r\'ealisation de Betti tangentielle} correspondante. Cette construction, men\'ee \`a bien dans \cite{Ay5}, repose sur un formidable appareil de techniques d\'evelopp\'ees ant\'erieurement par le m\^eme auteur, parmi lesquelles la d\'efinition des {\it cycles proches motiviques} et leur interpr\'etation dans le cadre d'une {\it variante rigide-analytique de} ${\bf{DM}}(\mathbb C((t)))$ \cite{Ay1}\cite{Ay4}. En adaptant les arguments du th. 4.5 ci-dessus \`a la r\'ealisation de Betti tangentielle et aux groupes de Galois motiviques {\og relatifs\fg} (via 5.1), il aboutit au terme d'un long parcours \`a l'analogue fonctionnel  du cor. 4.6 formul\'e au th. 0.1, \'enonc\'e inconditionnel d'une \'etonnante simplicit\'e\footnote{qui invite le lecteur optimiste \`a rechercher une preuve \'el\'ementaire, et le lecteur sceptique \`a tenter de fabriquer des contre-exemples - ce qui n'est d'ailleurs pas difficile si l'on s'autorise des fonctions alg\'ebriques ayant des singularit\'es sur le polydisque ({\og contre-exemples\fg} de m\^eme farine que ceux des ovales singuliers de Huygens et Leibniz, lecteurs sceptiques du lemme XXVIII de Newton).}.

  \smallskip
 {\it Remerciements.} Je remercie J. Ayoub, D. Bertrand et A. Huber pour leur aide durant la pr\'eparation de ce rapport, ainsi que P. Cartier pour sa relecture attentive.
 
 \end{sloppypar}

\end{document}